\documentclass{svproc}
\usepackage{amssymb}
\usepackage{amsmath}
\usepackage{amsfonts}
\usepackage{mathrsfs}
\newcommand{\R}{\mathbb{R}}
\newcommand{\sy}{\boldsymbol{\Psi}}
\newcommand{\py}{\boldsymbol{\Phi}}
\newtheorem{assumption}[proposition]{Assumption}
\newcommand{\inner}[2]{\langle #1, #2 \rangle}
\newcommand{\abs}[1]{\vert #1 \vert}
\newcommand{\norm}[1]{\Vert #1 \Vert}

\usepackage{url}

\begin{document}
\mainmatter              
\title{Existence and Uniqueness of Maximal Solutions to a 3D Navier-Stokes Equation with Stochastic Lie Transport}
\titlerunning{Solutions of the SALT Navier-Stokes Equation}  
%
\author{Daniel Goodair}
\authorrunning{Daniel Goodair} 
%
\tocauthor{Daniel Goodair}
\institute{Imperial College London, London SW7 2AZ, England,\\
\email{daniel.goodair16@imperial.ac.uk},\\ WWW home page:
\texttt{https://www.imperial.ac.uk/people/daniel.goodair16}}

\maketitle              

\begin{abstract}
We present here a criterion to conclude that an abstract SPDE posseses a unique maximal strong solution, which we apply to a three dimensional Stochastic Navier-Stokes Equation. Inspired by the work of Kato and Lai \cite{kato1984nonlinear} in the deterministic setting, we provide a comparable result here in the stochastic case whilst facilitating a variety of noise structures such as additive, multiplicative and transport. In particular our criterion is designed to fit viscous fluid dynamics models with Stochastic Advection by Lie Transport (SALT) as introduced in \cite{holm2015variational}. Our application to the Incompressible Navier-Stokes equation matches the existence and uniqueness result of the deterministic theory. This short work summarises the results and announces two papers \cite{Goodair1}, \cite{Goodair2} which give the full details for the abstract well-posedness arguments and application to the Navier-Stokes Equation respectively.

\keywords{stochastic transport, SPDE, Navier-Stokes, well-posedness}
\end{abstract}

\section{Introduction}
The theoretical analysis of fluid models perturbed by transport noise has been in significant demand since the release of the seminal works \cite{holm2015variational} and \cite{meminLU}. In the papers Holm and Memin establish a new class of stochastic equations driven by transport noise which serve as much improved fluid dynamics models by adding uncertainty in the transport of the fluid parcels to reflect the unresolved scales. Here we consider the SALT \cite{holm2015variational} Navier-Stokes Equation given by
\begin{equation} \label{number2equation}
    u_t - u_0 + \int_0^t\mathcal{L}_{u_s}u_s\,ds - \int_0^t \Delta u_s\, ds + \int_0^t Bu_s \circ d\mathcal{W}_s + \int_0^t \nabla \rho_s ds = 0
\end{equation}
and supplemented with the divergence-free (incompressibility) and zero-average conditions on the three dimensional torus $\mathbb{T}^3$. The equation is presented here in velocity form where $u$ represents the fluid velocity, $\rho$ the pressure, $\mathcal{L}$ is the mapping corresponding to the nonlinear term, $\mathcal{W}$ is a cylindrical Brownian Motion and $B$ is the relevant transport operator defined with respect to a collection of functions $(\xi_i)$ which physically represent spatial correlations. The explicit meaning of these conditions and the definitions of the operators involved are given at the beginning of Subsection \ref{ns sub}. These $(\xi_i)$ can be determined at coarse-grain resolutions from finely resolved numerical simulations, and mathematically are derived as eigenvectors of a velocity-velocity correlation matrix (see \cite{cotter1}, \cite{cotter2}, \cite{crisan3deuler}). The corresponding stochastic Euler equation was derived in \cite{streetvariational} and the viscous term plays no additional role in the stochastic derivation (without loss of generality we set the viscosity coefficient to be $1$).\\

 There has been limited progress in proving well-posedness for this class of equations: Crisan, Flandoli and Holm \cite{crisan3deuler} have shown local existence and uniqueness for the 3D Euler Equation on the torus, whilst Crisan and Lang (\cite{crisanrotating},\cite{crisanlang1},\cite{crisanlang2}) demonstrated the same result for the Euler, Rotating Shallow Water and Great Lake Equations on the torus once more. Whilst this represents a strong start in the theoretical analysis (alongside works for SPDEs with general transport noise e.g. \cite{Primitive}, \cite{general}), the modelling literature continues to expand in both the deterministic fluid models (see for example Figure 2 of \cite{crisanTQG} and the analysis therein) and method of stochastic perturbation (for example we may soon look to introduce nonlinearity and time dependence in the $(\xi_i)$). The significance of an abstract approach to the well-posedness question is clear, and whilst we discuss here only an application to SALT Navier-Stokes (\cite{holm2015variational},\cite{streetvariational}) the hope is that other stochastic viscous fluid models can be similarly solved by simply checking the required assumptions. We state our equation in the form
\begin{equation} \label{thespde}
    \sy_t = \sy_0 + \int_0^t \mathcal{A}(s,\sy_s)ds + \int_0^t\mathcal{G} (s,\sy_s) d\mathcal{W}_s
\end{equation}
for operators $\mathcal{A}$ and $\mathcal{G}$ to be elucidated in due course. The most notable contribution to the well-posedness theory for an abstract nonlinear SPDE is from \cite{temamabstract}. Here the authors prove the existence of a unique maximal solution to their abstract equation and apply this to the three dimensional primitive equations with a Lipschitz type multiplicative noise. The class of equations which we are concerned with include a differential operator in the noise term, preventing us from applying this framework. Moreover the assumptions on their operator $\mathcal{A}$ are quite explicit in terms of the sum of the standard fluid nonlinear term and a linear operator, which we don't restrict ourselves to. Overall our assumptions are much more general and allow for a straightforwards application to a wider class of SPDEs. Another relevant piece here is the work of Glatt-Holtz and Ziane \cite{holtzziane} whom show the same existence and uniqueness for the incompressible $3D$ Navier-Stokes with again a Lipschitz noise term. Though we cannot apply this method in the presence of our transport noise we look to adapt this argument to fit not just our Navier-Stokes equation but the wider class of stochastic viscous fluid models and SPDEs beyond. The impact of the boundary is fundamental to the equation and the approach of Glatt-Holtz and Ziane copes with the arising issues by working in the right function spaces; we recognised the importance of this in establishing an abstract framework which we hope to apply to such stochastic transport equations on the bounded domain as well.\\

This short summary work contains three more sections: in the subsequent one we properly define our Stochastic Navier-Stokes equation through the operators involved, the relevant function spaces, the notions of solution and main results. Following this we concretely define our abstract formulation and notion of solution, giving the assumptions that we require and the main results for the abstract equation. These assumptions are then all that needs to be checked to conclude the relevant existence and uniqueness for the proposed SPDE. In the final section we discuss the key steps behind proving these results; in the spirit of this as a summary work announcing our results we do not give a complete proof, though all such arguments are to be found in \cite{Goodair1}. We then address how our Navier-Stokes equation fits the context of the abstract formulation, though once more we do not give a thorough justification that the operators of our equation satisfy the required assumptions, with this precise treatment to come in \cite{Goodair2}.

\section{SALT Navier-Stokes and Results} \label{section salt navier stokes}

As alluded to in this section we formally introduce the equation (\ref{number2equation}) and state the main results.

\subsection{Preliminaries from Stochastic Analysis} \label{subsection prelims from stoch}

Throughout the paper we work with a fixed filtered probability space\\ $(\Omega,\mathcal{F},(\mathcal{F}_t), \mathbb{P})$ satisfying the usual conditions of completeness and right continuity. We take $\mathcal{W}$ to be a cylindrical Brownian Motion over some Hilbert Space $\mathfrak{U}$ with orthonormal basis $(e_i)$. The choice of $\mathfrak{U}$ and the subsequent basis play no role in the analysis. Recall (\cite{Goodair stoch}, Subsection 1.4) that $\mathcal{W}$ admits the representation $\mathcal{W}_t = \sum_{i=1}^\infty e_iW^i_t$ as a limit in $L^2(\Omega;\mathfrak{U}')$ whereby the $(W^i)$ are a collection of i.i.d. standard real valued Brownian Motions and $\mathfrak{U}'$ is an enlargement of the Hilbert Space $\mathfrak{U}$ such that the embedding $J: \mathfrak{U} \rightarrow \mathfrak{U}'$ is Hilbert-Schmidt and $\mathcal{W}$ is a $JJ^*-$cylindrical Brownian Motion over $\mathfrak{U}'$. Given a process $F:[0,T] \times \Omega \rightarrow \mathscr{L}^2(\mathfrak{U};\mathscr{H})$ progressively measurable and such that $F \in L^2\left(\Omega \times [0,T];\mathscr{L}^2(\mathfrak{U};\mathscr{H})\right)$, for any $0 \leq t \leq T$ we understand the stochastic integral $$\int_0^tF_sd\mathcal{W}_s$$ to be the infinite sum $$\sum_{i=1}^\infty \int_0^tF_s(e_i)dW^i_s$$ taken in $L^2(\Omega;\mathscr{H})$. We can extend this notion to processes $F$ which are such that $F(\omega) \in L^2\left( [0,T];\mathscr{L}^2(\mathfrak{U};\mathscr{H})\right)$ for $\mathbb{P}-a.e.$ $\omega$ via the traditional localisation procedure. In this case the stochastic integral is a local martingale in $\mathscr{H}$. A complete, direct construction of this integral, a treatment of its properties and the fundamentals of stochastic calculus in infinite dimensions can be found in \cite{Goodair stoch} Section 1. 

\subsection{SALT Navier-Stokes Equation} \label{ns sub}

We present the equation (\ref{number2equation}) on the three dimensional torus $\mathbb{T}^3$ (noting that all results hold on $\mathbb{T}^2$), and detail now the operators involved alongside the function spaces which define the equations. The operator $\mathcal{L}$ is defined for sufficiently regular functions $\phi,\psi:\mathbb{T}^3 \rightarrow \R^3$ by $$\mathcal{L}_{\phi}\psi:= \sum_{j=1}^3\phi^j\partial_j\psi$$ where $\phi^j:\mathbb{T}^3 \rightarrow \R$ is the $j^{\textnormal{th}}$ coordinate mapping of $\phi$ and $\partial_j\psi$ is defined by its $k^{\textnormal{th}}$ coordinate mapping $(\partial_j\psi)^k = \partial_j\psi^k$. The operator $B$ is defined as a linear operator on $\mathfrak{U}$ (introduced in $\ref{subsection prelims from stoch}$) by its action on the basis vectors $B(e_i,\cdot) := B_i(\cdot)$ by $$B_i = \mathcal{L}_{\xi_i} + \mathcal{T}_{\xi_i}$$ for $\mathcal{L}$ as above and $$\mathcal{T}_{\phi}\psi:= \sum_{j=1}^3\psi^j\nabla\phi^j.$$  A complete discussion of how $B$ is then defined on $\mathfrak{U}$ is given in \cite{Goodair stoch} Subsection 2.2. We embed the divergence-free and zero-average conditions into the relevant function spaces and simply define our solutions as belonging to these spaces. To be explicit, by a divergence-free function we mean a $\phi \in W^{1,2}(\mathbb{T}^3;\R^3)$ such that $$\sum_{j=1}^3\partial_j\phi^j = 0$$ and by zero-average we ask for a $\psi \in L^2(\mathbb{T}^3;\R^3)$ with the property $$\int_{\mathbb{T}^3}\psi \ d\lambda = 0$$ for $\lambda$ the Lebesgue measure on $\mathbb{T}^3$. We introduce the space $L^2_{\sigma}(\mathbb{T}^3;\R^3)$ as the subspace of $L^2(\mathbb{T}^3;\R^3)$ consisting of zero-average functions which are 'weakly divergence-free'; see \cite{robinson} Definition 2.1 for the precise construction. $W^{1,2}_{\sigma}(\mathbb{T}^3;\R^3)$ is then defined as the subspace of $W^{1,2}(\mathbb{T}^3;\R^3)$ consisting of zero-average\\ divergence-free functions, and $W^{2,2}_{\sigma}(\mathbb{T}^3;\R^3):=W^{2,2}(\mathbb{T}^3;\R^3) \cap W^{1,2}_{\sigma}(\mathbb{T}^3;\R^3)$.\\

As is standard in the treatment of the incompressible Navier-Stokes Equation we consider a projected version to eliminate the pressure term and facilitate us working in the above spaces. Note that $\rho$ does not come with an evolution equation and is simply chosen to ensure the incompressibility condition. The idea is to solve the projected equation and then append a pressure to it, see \cite{robinson}. To this end we introduce the standard Leray Projector $\mathcal{P}$ defined as the orthogonal projection in $L^2(\mathbb{T}^3;\R^3)$ onto $L^2_{\sigma}(\mathbb{T}^3;\R^3)$. As we look to project the equation (\ref{number2equation}) as discussed, we ought to address the Stratonovich integral. We look to convert this term into an It\^{o} integral to enable our analysis, but the resulting converted and projected equation should not depend on the order in which the projection and conversion occur. To this end we assume that the $(\xi_i)$ are such that $\xi_i \in W^{1,2}_{\sigma}(\mathbb{T}^3;\R^3) \cap W^{3,\infty}(\mathbb{T}^3;\R^3)$ and satisfy the bound \begin{equation} \label{xi bound}
    \sum_{i=1}^\infty\norm{\xi_i}_{W^{3,\infty}}^2 < \infty.
\end{equation}
The significance of the bound (\ref{xi bound}) will be revisited, but for now we note that as each $\xi_i$ is divergence-free then each $B_i$ satisfies the property that $\mathcal{P}B_i$ is equal to $\mathcal{P}B_i\mathcal{P}$ on $W^{1,2}(\mathbb{T}^3;\R^3)$ which ensures that the projection and conversion commute. Our new equation is then
\begin{align} \nonumber
     u_t - u_0 + \int_0^t\mathcal{P}\mathcal{L}_{u_s}u_s\,ds &+ \int_0^t A u_sds \\&-  \frac{1}{2}\sum_{i=1}^\infty \int_0^t\mathcal{P}B_i^2u_sds + \sum_{i=1}^\infty\int_0^t \mathcal{P}B_iu_s  dW_s^i = 0 \label{converted equation}
\end{align}
where $A := -\mathcal{P}\Delta$ is known as the Stokes Operator. Details of the It\^{o}-Stratonovich conversion can be found in \cite{Goodair stoch} Subsection 2.3. We shall use the Stokes operator to define inner products with which we equip our function spaces. Recall from \cite{robinson} Theorem 2.24 for example that there exists a collection of functions $(a_k)$, $a_k \in W^{1,2}_{\sigma}(\mathbb{T}^3;\R^3) \cap C^{\infty}(\mathbb{T}^3;\R^3)$ such that the $(a_k)$ are eigenfunctions of $A$, are an orthonormal basis in $L^2_{\sigma}(\mathbb{T}^3;\R^3)$ and an orthogonal basis in $W^{1,2}_{\sigma}(\mathbb{T}^3;\R^3)$ considered as Hilbert Spaces with standard $L^2(\mathbb{T}^3;\R^3)$, $W^{1,2}(\mathbb{T}^3;\R^3)$ inner products. The corresponding eigenvalues $(\lambda_k)$ are strictly positive and approach infinity as $k \rightarrow \infty$. Thus any $\phi \in W^{1,2}_{\sigma}(\mathbb{T}^3;\R^3)$ admits the representation $$\phi = \sum_{k=1}^\infty \phi_ka_k$$ so for $m \in \mathbb{N}$ we can define $A^{m/2}$ by $$A^{m/2}: \phi \mapsto  \sum_{k=1}^\infty \lambda_k^{m/2}\phi_ka_k$$ which is a well defined element of $L^2_{\sigma}(\mathbb{T}^3;\R^3)$ on any $\phi$
such that \begin{equation} \label{powerproperty} \sum_{k=1}^\infty \lambda_k^m\phi_k^2 < \infty.\end{equation} For $\phi,\psi$ with the property (\ref{powerproperty}) then the bilinear form $$\inner{\phi}{\psi}_{m}:= \inner{A^{m/2}\phi}{A^{m/2}\psi}$$ is well defined. For $m=1,2$ this is an inner product on the spaces $W^{1,2}_{\sigma}(\mathbb{T}^3;\R^3)$, $W^{2,2}_{\sigma}(\mathbb{T}^3;\R^3)$ respectively which is equivalent to the standard $W^{1,2}(\mathbb{T}^3;\R^3)$, $W^{2,2}(\mathbb{T}^3;\R^3)$ inner product. Of course $\inner{\cdot}{\cdot}_3$ is well defined on\\ $\bigcup_{k=1}^\infty\textnormal{span}\{a_1,\dots,a_k\}$ and so we define $W^{3,2}_{\sigma}(\mathbb{T}^3;\R^3)$ as the completion of\\ $\bigcup_{k=1}^\infty\textnormal{span}\{a_1,\dots,a_k\}$ in this inner product. We consider $W^{m,2}_{\sigma}(\mathbb{T}^3;\R^3)$ as a Hilbert Space equipped with the $\inner{\cdot}{\cdot}_m$ inner product, and define our solution to the equation (\ref{converted equation}) relative to these spaces.

\subsection{Notions of Solution and Results}

We frame this definition for an $\mathcal{F}_0-$measurable $u_0: \Omega \rightarrow W^{1,2}_{\sigma}(\mathbb{T}^3;\R^3)$. Here and throughout we use the notation $\mathbf{1}$ for the indicator function.

\begin{definition} \label{definitionofNS solution}
A pair $(u,\tau)$ where $\tau$ is a $\mathbb{P}-a.s.$ positive stopping time and $u$ is a process such that for $\mathbb{P}-a.e.$ $\omega$, $u_{\cdot}(\omega) \in C\left([0,T];W^{1,2}_{\sigma}(\mathbb{T}^3;\R^3)\right)$ and $u_{\cdot}(\omega)\mathbf{1}_{\cdot \leq \tau(\omega)} \in L^2\left([0,T];W^{2,2}_{\sigma}(\mathbb{T}^3;\R^3)\right)$ for all $T>0$ with $u_{\cdot}\mathbf{1}_{\cdot \leq \tau}$ progressively measurable in $W^{2,2}_{\sigma}(\mathbb{T}^3;\R^3)$, is said to be a local strong solution of the equation (\ref{thespde}) if the identity
\begin{align} \nonumber
     u_t - u_0 &+ \int_0^{t\wedge \tau}\mathcal{P}\mathcal{L}_{u_s}u_s\,ds + \int_0^{t\wedge \tau} A u_sds \\& \qquad -  \frac{1}{2}\sum_{i=1}^\infty \int_0^{t\wedge \tau}\mathcal{P}B_i^2u_sds + \sum_{i=1}^\infty\int_0^{t\wedge \tau} \mathcal{P}B_iu_s  dW_s^i = 0 \label{converted equation identity}
\end{align}
holds $\mathbb{P}-a.s.$ in $L^2_{\sigma}(\mathbb{T}^3;\R^3)$ for all $t \geq 0$.
\end{definition}
We shall address why this definition makes sense in the abstract setting in Subsection \ref{subsection abstract notions of solution and results}, before then translating this abstract framework back to our Navier-Stokes Equation.  
\begin{definition} \label{NS maximal definition}
A pair $(u,\Theta)$ such that there exists a sequence of stopping times $(\theta_j)$ which are $\mathbb{P}-a.s.$ monotone increasing and convergent to $\Theta$, whereby $(u_{\cdot \wedge \theta_j},\theta_j)$ is a local strong solution of the equation (\ref{converted equation}) for each $j$, is said to be a maximal strong solution of the equation (\ref{converted equation}) if for any other pair $(v,\Gamma)$ with this property then $\Theta \leq \Gamma$ $\mathbb{P}-a.s.$ implies $\Theta = \Gamma$ $\mathbb{P}-a.s.$.
\end{definition}

\begin{definition} \label{NS unique}
A maximal strong solution $(u,\Theta)$ of the equation (\ref{converted equation}) is said to be unique if for any other such solution $(v,\Gamma)$, then $\Theta = \Gamma$ $\mathbb{P}-a.s.$ and for all $t \in [0,\Theta)$, \begin{equation} \nonumber\mathbb{P}\left(\left\{\omega \in \Omega: u_{t}(\omega) =  v_{t}(\omega)  \right\} \right) = 1. \end{equation}
\end{definition}
We can now state the main result of the paper.

\begin{theorem} \label{NS theorem}
For any given $\mathcal{F}_0-$ measurable $u_0:\Omega \rightarrow W^{1,2}_{\sigma}(\mathbb{T}^3;\R^3)$, there exists a unique maximal strong solution $(u,\Theta)$ of the equation (\ref{converted equation}). Moreover at $\mathbb{P}-a.e.$ $\omega$ for which $\Theta(\omega)<\infty$, we have that \begin{equation}\label{blowuppropertyNS}\sup_{r \in [0,\Theta(\omega))}\norm{u_r(\omega)}_1^2 + \int_0^{\Theta(\omega)}\norm{u_r(\omega)}_2^2dr = \infty.\end{equation}
\end{theorem}

\section{Abstract Framework and Results}

We now establish the abstract framework through which we arrive at Theorem \ref{NS theorem}. This involves giving two sets of assumptions before exploring the abstract method with the assumptions in place, and then in Section \ref{section salt navier stokes in abstract framework} discussing how (\ref{converted equation}) fits into this framework. These assumption sets pertain to two different notions of solution (both strong in the probabilistic sense but related to different spaces), the reason for which will be illustrated in Section \ref{Section abstract solution method}. We give these as two distinct sets of assumptions in the event that an equation fits the first set of assumptions but not the second, such that we would still be able to conclude that some type of solution exists for the equation.  
\subsection{Assumption Set 1} \label{subsection assumption set 1}

We work with a quartet of continuously embedded Hilbert Spaces $$V \hookrightarrow H \hookrightarrow U \hookrightarrow X$$ and the operators \begin{align*}
    \mathcal{A}&:[0,\infty) \times V \rightarrow U,\\
    \mathcal{G}&:[0,\infty) \times V \rightarrow \mathscr{L}^2(\mathfrak{U};H).
\end{align*} 
We ask that there is a continuous bilinear form $\inner{\cdot}{\cdot}_{X \times H}: X \times H \rightarrow \R$ such that for $\phi \in U$ and $\psi \in H$, \begin{equation} \label{bilinear form}
    \inner{\phi}{\psi}_{X \times H} =  \inner{\phi}{\psi}_{U}.
\end{equation}
Moreover the continuity and bilinearity ensures that there exists some constant $c$ whereby for all such $\phi,\psi$, \begin{equation}
    \label{cauchy schwartz for bilinear form}
    \abs{\inner{\phi}{\psi}_{X \times H}} \leq c\norm{\phi}_X{\norm{\psi}}_H.
\end{equation}
As we look to use a Galerkin Scheme to solve our equation, we introduce now a sequence of spaces $(V_n)$ contained in $V$ given by $V_n:= \textnormal{span}\left\{a_1, \dots, a_n \right\}$ for $(a_n)$ an orthogonal basis in $U$. Defining $\mathcal{P}_n$ to be the orthogonal projection onto $V_n$ in $X$, we shall also assume that the restriction of $\mathcal{P}_n$ to $U$ is an orthogonal projection in $U$ and that the sequence of these projections is uniformly bounded on $H$: that is, that there exists some constant $c$ independent of $n$ such that for all $\boldsymbol{\phi} \in H$,

\begin{equation} \label{projectionsboundedonH}
    \norm{\mathcal{P}_n\boldsymbol{\phi}}_H^2 \leq c\norm{\boldsymbol{\phi}}_H^2.
\end{equation}
We also require the existence of a real valued sequence $(\mu_n)$ with $\mu_n \rightarrow \infty$, which is such that for any $\phi \in U$ and $\psi \in H$, \begin{align}
    \label{mu1}\norm{(I - \mathcal{P}_n)\phi}_X \leq \frac{1}{\mu_n}\norm{\phi}_U,\\ \label{mu2}
    \norm{(I - \mathcal{P}_n)\psi}_U \leq \frac{1}{\mu_n}\norm{\psi}_H
\end{align}
where $I$ represents the identity operator in the corresponding spaces. These assumptions are of course supplemented by a series of assumptions on the operators. We shall use general notation $c_t$ to represent a function $c_\cdot:[0,\infty) \rightarrow \R$ bounded on $[0,T]$ for any $T > 0$, evaluated at the time $t$. Moreover we define functions $K$, $\tilde{K}$ relative to some non-negative constants $p,\tilde{p},q,\tilde{q}$. We use a generic notation to define the functions $K: U \rightarrow \R$, $K: U \times U \rightarrow \R$, $\tilde{K}: H \rightarrow \R$ and $\tilde{K}: H \times H \rightarrow \R$ by
\begin{align*}
    K(\phi)&:= 1 + \norm{\phi}_U^{p},\\
    K(\phi,\psi)&:= 1+\norm{\phi}_U^{p} + \norm{\psi}_U^{q},\\
    \tilde{K}(\phi) &:= K(\phi) + \norm{\phi}_H^{\tilde{p}},\\
    \tilde{K}(\phi,\psi) &:= K(\phi,\psi) + \norm{\phi}_H^{\tilde{p}} + \norm{\psi}_H^{\tilde{q}}
\end{align*}
 Distinct use of the function $K$ will depend on different constants but in no meaningful way in our applications, hence no explicit reference to them shall be made. In the case of $\tilde{K}$, when $\tilde{p}, \tilde{q} = 2$ then we shall denote the general $\tilde{K}$ by $\tilde{K}_2$. In this case no further assumptions are made on the $p,q$. That is, $\tilde{K}_2$ has the general representation \begin{equation}\label{Ktilde2}\tilde{K}_2(\phi,\psi) = K(\phi,\psi) + \norm{\phi}_H^2 + \norm{\psi}_H^2\end{equation} and similarly as a function of one variable.\\
 
 We state the assumptions for arbitrary elements $\phi,\psi \in V$, $\phi^n \in V_n$ and $t \in [0,\infty)$, and a fixed $\kappa > 0$. Understanding $\mathcal{G}$ as an operator $\mathcal{G}: [0,\infty) \times V \times \mathfrak{U} \rightarrow H$, we introduce the notation $\mathcal{G}_i(\cdot,\cdot):= \mathcal{G}(\cdot,\cdot,e_i)$.
 
 \begin{assumption} \label{measurabilityassumption}
For any $T>0$, $\mathcal{A}:[0,T] \times V \rightarrow U$ and  $\mathcal{G}:[0,T] \times V \rightarrow \mathscr{L}^2(\mathfrak{U};H)$ are measurable. 
\end{assumption}

\begin{remark}
Measurability here and throughout the paper is defined with respect to the Borel Sigma Algebra on the relevant Hilbert Spaces.
\end{remark}
 
  \begin{assumption} \label{new assumption 1}
 \begin{align}
     \label {111} \norm{\mathcal{A}(t,\boldsymbol{\phi})}^2_U +\sum_{i=1}^\infty \norm{\mathcal{G}_i(t,\boldsymbol{\phi})}^2_H &\leq c_t K(\boldsymbol{\phi})\left[1 + \norm{\boldsymbol{\phi}}_V^2\right],\\ \label{222}
     \norm{\mathcal{A}(t,\boldsymbol{\phi}) - \mathcal{A}(t,\boldsymbol{\psi})}_X &\leq  c_t\left[K(\phi,\psi) + \norm{\phi}_V + \norm{\psi}_V \right]\norm{\phi-\psi}_H,\\ \label{333}
    \sum_{i=1}^\infty \norm{\mathcal{G}_i(t,\boldsymbol{\phi}) - \mathcal{G}_i(t,\boldsymbol{\psi})}_X &\leq c_tK(\phi,\psi)\norm{\phi-\psi}_H
 \end{align}
 \end{assumption}

\begin{assumption} \label{assumptions for uniform bounds2}
\begin{align}
   \label{uniformboundsassumpt1} \nonumber 2\inner{\mathcal{P}_n\mathcal{A}(t,\boldsymbol{\phi}^n)}{\boldsymbol{\phi}^n}_H + \sum_{i=1}^\infty\norm{\mathcal{P}_n\mathcal{G}_i(t,\boldsymbol{\phi}^n)}_H^2 &\leq\\ c_t\tilde{K}_2(\boldsymbol{\phi}^n)&\left[1 + \norm{\boldsymbol{\phi}^n}_H^2\right] - \kappa\norm{\boldsymbol{\phi}^n}_V^2,\\  \label{uniformboundsassumpt2}
    \sum_{i=1}^\infty \inner{\mathcal{P}_n\mathcal{G}_i(t,\boldsymbol{\phi}^n)}{\boldsymbol{\phi}^n}^2_H &\leq c_t\tilde{K}_2(\boldsymbol{\phi}^n)\left[1 + \norm{\boldsymbol{\phi}^n}_H^4\right].
\end{align}
\end{assumption}

\begin{assumption} \label{therealcauchy assumptions}
\begin{align}
  \nonumber 2\inner{\mathcal{A}(t,\boldsymbol{\phi}) - \mathcal{A}(t,\boldsymbol{\psi})}{\boldsymbol{\phi} - \boldsymbol{\psi}}_U &+ \sum_{i=1}^\infty\norm{\mathcal{G}_i(t,\boldsymbol{\phi}) - \mathcal{G}_i(t,\boldsymbol{\psi})}_U^2\\ \label{therealcauchy1} &\leq  c_{t}\tilde{K}_2(\boldsymbol{\phi},\boldsymbol{\psi}) \norm{\boldsymbol{\phi}-\boldsymbol{\psi}}_U^2 - \kappa\norm{\boldsymbol{\phi}-\boldsymbol{\psi}}_H^2,\\ \label{therealcauchy2}
    \sum_{i=1}^\infty \inner{\mathcal{G}_i(t,\boldsymbol{\phi}) - \mathcal{G}_i(t,\boldsymbol{\psi})}{\boldsymbol{\phi}-\boldsymbol{\psi}}^2_U & \leq c_{t} \tilde{K}_2(\boldsymbol{\phi},\boldsymbol{\psi}) \norm{\boldsymbol{\phi}-\boldsymbol{\psi}}_U^4
\end{align}
\end{assumption}

\begin{assumption} \label{assumption for prob in V}
\begin{align}
   \label{probability first} 2\inner{\mathcal{A}(t,\boldsymbol{\phi})}{\boldsymbol{\phi}}_U + \sum_{i=1}^\infty\norm{\mathcal{G}_i(t,\boldsymbol{\phi})}_U^2 &\leq c_tK(\boldsymbol{\phi})\left[1 +  \norm{\boldsymbol{\phi}}_H^2\right],\\\label{probability second}
    \sum_{i=1}^\infty \inner{\mathcal{G}_i(t,\boldsymbol{\phi})}{\boldsymbol{\phi}}^2_U &\leq c_tK(\boldsymbol{\phi})\left[1 + \norm{\boldsymbol{\phi}}_H^4\right].
\end{align}
\end{assumption}

\subsection{Assumption Set 2}
 
These assumptions are only checked in addition to Assumption Set 1 and so take place in the same framework. We state the assumptions now for arbitrary elements $\phi,\psi \in H$ and $t \in [0,\infty)$, and continue to use the $c,K,\tilde{K}, \kappa$ notation of Assumption Set 1. 

\begin{assumption} \label{secondmeasurabilityassumption}
For any $T>0$, $\mathcal{A}:[0,T] \times H \rightarrow X$ is measurable, and whenever $\py$ is a progressively measurable process in $H$ we have that $\mathcal{G}(\cdot,\py_{\cdot})$ is progressively measurable in $\mathscr{L}^2(\mathfrak{U};U)$. 
\end{assumption}

\begin{assumption}
\begin{align}
    \label{wellposedinX}
    \norm{\mathcal{A}(t,\boldsymbol{\phi})}^2_X +  \sum_{i=1}^\infty \norm{\mathcal{G}_i(t,\boldsymbol{\phi})}^2_U &\leq c_tK(\boldsymbol{\phi})\left[1 + \norm{\boldsymbol{\phi}}_H^2\right],\\ \label{222*} \norm{\mathcal{A}(t,\boldsymbol{\phi}) - \mathcal{A}(t,\boldsymbol{\psi})}_X &\leq  c_t\left[K(\phi,\psi) + \norm{\phi}_H + \norm{\psi}_H \right]\norm{\phi-\psi}_H
\end{align}

\end{assumption}

\begin{assumption} \label{uniqueness for H valued}
\begin{align}
    \nonumber 2\inner{\mathcal{A}(t,\boldsymbol{\phi}) - \mathcal{A}(t,\boldsymbol{\psi})}{\boldsymbol{\phi} - \boldsymbol{\psi}}_X + \sum_{i=1}^\infty\norm{\mathcal{G}_i(t,\boldsymbol{\phi}) - \mathcal{G}_i(t,\boldsymbol{\psi})}_X^2 &\leq\\ \label{therealcauchy1*} c_{t}\tilde{K}_2(\boldsymbol{\phi},\boldsymbol{\psi}) \norm{\boldsymbol{\phi}-\boldsymbol{\psi}}_X^2,&\\ \nonumber
    \sum_{i=1}^\infty \inner{\mathcal{G}_i(t,\boldsymbol{\phi}) - \mathcal{G}_i(t,\boldsymbol{\psi})}{\boldsymbol{\phi}-\boldsymbol{\psi}}^2_X & \leq\\\label{therealcauchy2*} c_{t} \tilde{K}_2(\boldsymbol{\phi},\boldsymbol{\psi}) \norm{\boldsymbol{\phi}-\boldsymbol{\psi}}_X^4&
\end{align}
\end{assumption}
We in fact state Assumption \ref{assumption for probability in H} for $\phi \in V$ and some $\kappa > 0$, making this a stronger assumption than \ref{assumption for prob in V}.

\begin{assumption} \label{assumption for probability in H}
With the stricter requirement that $\phi\in V$ then 
\begin{align}
   \label{probability first H} 2\inner{\mathcal{A}(t,\boldsymbol{\phi})}{\boldsymbol{\phi}}_U + \sum_{i=1}^\infty\norm{\mathcal{G}_i(t,\boldsymbol{\phi})}_U^2 &\leq c_tK(\boldsymbol{\phi}) -  \kappa\norm{\boldsymbol{\phi}}_H^2,\\\label{probability second H}
    \sum_{i=1}^\infty \inner{\mathcal{G}_i(t,\boldsymbol{\phi})}{\boldsymbol{\phi}}^2_U &\leq c_tK(\boldsymbol{\phi}).
\end{align}
\end{assumption}

\subsection{Notions of Solution and Results} \label{subsection abstract notions of solution and results}

Here we define the two different notions of solution, which we call $V-$valued solutions and $H-$valued solutions. The corresponding definitions of uniqueness and maximality are given in one for both notions of solution. We frame the definition of the $V-$valued solutions for an initial condition $\sy_0:\Omega \rightarrow H$ which is an $\mathcal{F}_0-$measurable mapping, and for the $H-$valued solutions a $\sy_0:\Omega \rightarrow U$ which is likewise $\mathcal{F}_0-$measurable.

\begin{definition} \label{definitionofregularsolution}
A pair $(\sy,\tau)$ where $\tau$ is a $\mathbb{P}-a.s.$ positive stopping time and $\sy$ is a process such that for $\mathbb{P}-a.e.$ $\omega$, $\sy_{\cdot}(\omega) \in C\left([0,T];H\right)$ and $\sy_{\cdot}(\omega)\mathbf{1}_{\cdot \leq \tau(\omega)} \in L^2\left([0,T];V\right)$ for all $T>0$ with $\sy_{\cdot}\mathbf{1}_{\cdot \leq \tau}$ progressively measurable in $V$, is said to be a $V$-valued local strong solution of the equation (\ref{thespde}) if the identity
\begin{equation} \label{identityindefinitionoflocalsolution}
    \sy_{t} = \sy_0 + \int_0^{t\wedge \tau} \mathcal{A}(s,\sy_s)ds + \int_0^{t \wedge \tau}\mathcal{G} (s,\sy_s) d\mathcal{W}_s
\end{equation}
holds $\mathbb{P}-a.s.$ in $U$ for all $t \geq 0$.
\end{definition}

\begin{remark} \label{remark1}
If $(\sy,\tau)$ is a $V$-valued local strong solution of the equation (\ref{thespde}), then $\sy_\cdot = \sy_{
\cdot \wedge \tau}$.
\end{remark}

\begin{remark} \label{remark on prog meas equivalence}
The progressive measurability condition on $\sy_{\cdot}\mathbf{1}_{\cdot \leq \tau}$ may look a little suspect as $\sy_0$ itself may only belong to $H$ and not $V$ making it impossible for $\sy_{\cdot}\mathbf{1}_{\cdot \leq \tau}$ to be even adapted in $V$. We are mildly abusing notation here; what we really ask is that there exists a process $\py$ which is progressively measurable in $V$ and such that $\py_{\cdot} = \sy_{\cdot}\mathbf{1}_{\cdot \leq \tau}$ almost surely over the product space $\Omega \times [0,\infty)$ with product measure $\mathbb{P}\times \lambda$ for $\lambda$ the Lebesgue measure on $[0,\infty)$. 
\end{remark}

\begin{remark}
If Assumption \ref{measurabilityassumption} and (\ref{111}) hold, then  the time integral is well defined in $U$ and the stochastic integral is well defined as a local martingale in $H$.
\end{remark}
 
 \begin{definition} \label{definitionofHsolution}
A pair $(\sy,\tau)$ where $\tau$ is a $\mathbb{P}-a.s.$ positive stopping time and $\sy$ is a process such that for $\mathbb{P}-a.e.$ $\omega$, $\sy_{\cdot}(\omega) \in C\left([0,T];U\right)$ and $\sy_{\cdot}(\omega)\mathbf{1}_{\cdot \leq \tau(\omega)} \in L^2\left([0,T];H\right)$ for all $T>0$ with $\sy_{\cdot}\mathbf{1}_{\cdot \leq \tau}$ progressively measurable in $H$, is said to be an $H$-valued local strong solution of the equation (\ref{thespde}) if the identity
\begin{equation} \label{identityindefinitionoflocalsolution**}
    \sy_{t} = \sy_0 + \int_0^{t\wedge \tau} \mathcal{A}(s,\sy_s)ds + \int_0^{t \wedge \tau}\mathcal{G} (s,\sy_s) d\mathcal{W}_s
\end{equation}
holds $\mathbb{P}-a.s.$ in $X$ for all $t \geq 0$.
\end{definition}
 
\begin{remark} \label{remark analogy}
The analogy to Remarks \ref{remark1}, \ref{remark on prog meas equivalence} hold for the $H-$valued solutions.
\end{remark}

\begin{remark}
If Assumption \ref{secondmeasurabilityassumption} and (\ref{wellposedinX}) hold, then  the time integral is well defined in $X$ and the stochastic integral is well defined as a local martingale in $U$.
\end{remark}

In the following we use $V;H$ to mean $V$ or $H$ respectively.

\begin{definition} \label{V valued maximal definition}
A pair $(\sy,\Theta)$ such that there exists a sequence of stopping times $(\theta_j)$ which are $\mathbb{P}-a.s.$ monotone increasing and convergent to $\Theta$, whereby $(\sy_{\cdot \wedge \theta_j},\theta_j)$ is a $(V;H)-$valued local strong solution of the equation (\ref{thespde}) for each $j$, is said to be a $(V;H)-$valued maximal strong solution of the equation (\ref{thespde}) if for any other pair $(\py,\Gamma)$ with this property then $\Theta \leq \Gamma$ $\mathbb{P}-a.s.$ implies $\Theta = \Gamma$ $\mathbb{P}-a.s.$.
\end{definition}

\begin{definition}
A $(V;H)-$valued maximal strong solution $(\sy,\Theta)$ of the equation (\ref{thespde}) is said to be unique if for any other such solution $(\py,\Gamma)$, then $\Theta = \Gamma$ $\mathbb{P}-a.s.$ and for all $t \in [0,\Theta)$, \begin{equation} \nonumber\mathbb{P}\left(\left\{\omega \in \Omega: \sy_{t}(\omega) =  \py_{t}(\omega)  \right\} \right) = 1. \end{equation}
\end{definition}

\begin{theorem} \label{theorem1}
Suppose that Assumption Set 1 holds. Then for any given $\mathcal{F}_0-$ measurable $\sy_0:\Omega \rightarrow H$, there exists a unique $V-$valued maximal strong solution $(\sy,\Theta)$ of the equation (\ref{thespde}). Moreover at $\mathbb{P}-a.e.$ $\omega$ for which $\Theta(\omega)<\infty$, we have that \begin{equation}\label{blowupproperty}\sup_{r \in [0,\Theta(\omega))}\norm{\sy_r(\omega)}_H^2 + \int_0^{\Theta(\omega)}\norm{\sy_r(\omega)}_V^2dr = \infty.\end{equation}
\end{theorem}

\begin{theorem} \label{theorem2}
Suppose that Assumption Set 1 and 2 hold. Then for any given $\mathcal{F}_0-$measurable $\sy_0:\Omega \rightarrow U$, there exists a unique $H-$valued maximal strong solution $(\sy,\Theta)$ of the equation (\ref{thespde}). Moreover at $\mathbb{P}-a.e.$ $\omega$ for which $\Theta(\omega)<\infty$, we have that \begin{equation}\label{blowuppropertyH}\sup_{r \in [0,\Theta(\omega))}\norm{\sy_r(\omega)}_U^2 + \int_0^{\Theta(\omega)}\norm{\sy_r(\omega)}_H^2dr = \infty.\end{equation}
\end{theorem}

\section{Abstract Solution Method and Application} \label{Section abstract solution method}

In this final section we give the main steps of the proofs of Theorems \ref{theorem1} and \ref{theorem2}, followed by a brief exposition of how our SALT Navier-Stokes Equation fits into this framework.

\subsection{Abstract Solution Method}

\begin{proof}[Theorem \ref{theorem1}] We suppose that Assumption Set 1 holds and address the question first for an initial condition $\sy_0$ which is such that for $\mathbb{P}-a.e.$ $\omega$, \begin{equation} \label{boundedinitialcondition} \norm{\sy_0(\omega)}^2_H \leq M'\end{equation} for some constant $M'$. We work with this bounded initial condition in the first instance as we shall use local solutions up to first hitting times given in terms of the initial condition, so this boundedness translates to boundedness of the relevant process up until these times. As directed in Subsection \ref{subsection assumption set 1} we are to use a Galerkin Scheme, whereby we consider the equations \begin{equation} \label{nthorderGalerkin}
       \sy^n_t = \sy^n_0 + \int_0^t \mathcal{P}_n\mathcal{A}(s,\sy^n_s)ds + \int_0^t\mathcal{P}_n\mathcal{G} (s,\sy^n_s) d\mathcal{W}_s
\end{equation}
with notation $\mathcal{P}_n\mathcal{G} (\cdot,\cdot,e_i):=\mathcal{P}_n\mathcal{G}_{i}(\cdot,\cdot).$ A local strong solution of this equation is defined as a  pair $(\sy^n,\tau)$ where $\tau$ is a $\mathbb{P}-a.s.$ positive stopping time and $\sy^n$ is an adapted process in $V_n$ such that for $\mathbb{P}-a.e.$ $\omega$, $\sy^n_{\cdot}(\omega) \in C\left([0,T];V_n\right)$ for all $T>0$, and the identity
\begin{equation} \label{identityindefinitionoflocalgalerkinsolution}
    \sy^n_{t} = \sy^n_0 + \int_0^{t\wedge \tau} \mathcal{P}_n\mathcal{A}(s,\sy^n_s)ds + \int_0^{t \wedge \tau}\mathcal{P}_n\mathcal{G} (s,\sy^n_s) d\mathcal{W}_s
\end{equation}
holds $\mathbb{P}-a.s.$ in $V_n$ for all $t \geq 0$. We can conclude that for any fixed $t > 0$ and $M>1$, a local strong solution $(\sy^n,\tau^{M,t}_n)$ of (\ref{nthorderGalerkin}) exists for the stopping time $\tau^{M,t}_n$ defined by \begin{equation}\label{tauMtn}\tau^{M,t}_n := t \wedge \inf\left\{s \geq 0: \sup_{r \in [0,s]}\norm{\sy^{n}_{r}}^2_U + \int_0^s\norm{\sy^{n}_{r}}^2_Hdr \geq M + \norm{\sy^n_0}_U^2 \right\}.\end{equation}
This conclusion is reached thanks to Assumption \ref{new assumption 1}, through standard theory in the finite dimensional Hilbert Space $V_n$ though some care must be taken for the infinite dimensional Brownian Motion. Understanding that \begin{equation} \label{uniformboundsofinitialconditioninH}\norm{\sy^n_0(\omega)}^2_H \leq c\norm{\sy_0(\omega)}^2_H \leq cM'\end{equation} coming from (\ref{projectionsboundedonH}) and (\ref{boundedinitialcondition}), it is clear that \begin{equation} \label{uniformboundofinitialcondition}
    \norm{\sy^n_0(\omega)}^2_U \leq \tilde{M}
    \end{equation}
for some $\tilde{M}$ clearly still independent of $n$ and $\omega$. Thus we see the bound \begin{equation} \label{galerkinboundstoppingtime}
    \sup_{r \in [0,\tau^{M,t}_n(\omega)]}\norm{\sy^n_{r}(\omega)}_U^2 + \int_0^{\tau^{M,t}_n(\omega)}\norm{\sy^n_s(\omega)}_H^2ds \leq M + \tilde{M}
\end{equation} holds true for every $n$ and $\mathbb{P}-a.e.$ $\omega$. This boundedness plays a significant role in our analysis and demonstrates the importance of starting from this bounded initial condition in the first instance. The motivation for choosing these stopping times comes from the work of Glatt-Holtz and Ziane in the referenced paper \cite{holtzziane}. The authors prove an abstract result which is the central theorem of the paper, which we simply restate in the Appendix as Theorem \ref{greenlemma}. In the original paper, the authors use the traditional Galerkin Scheme for Navier-Stokes (given by the basis of eigenfunctions of the Stokes Operator) and apply this theorem directly with the spaces $\mathcal{H}_1:= W^{2,2}_{\sigma}(\mathbb{T}^3;\R^3)$, $\mathcal{H}_2:= W^{1,2}_{\sigma}(\mathbb{T}^3;\R^3)$. We have to take a slight detour from this method in the case of transport noise due to the condition (\ref{required cauchy}). Translating this to our framework through $\mathcal{H}_1 = H$ and $\mathcal{H}_2 = U$, the idea in showing this condition is to apply the It\^{o} Formula in $U$ to the difference process $\sy^n-\sy^m$. When we simplify down the term arising from the quadratic variation of the stochastic integral, we must control $$\sum_{i=1}^\infty\norm{[I-\mathcal{P}_m]\mathcal{G}_i(s,\sy^m_s)}_{U}^2$$  which we would do via (\ref{mu2}) and (\ref{projectionsboundedonH}) to bound the above by $$\sum_{i=1}^\infty\frac{1}{\mu_m}\norm{\mathcal{G}_i(s,\sy^m_s)}_{H}^2.$$ In order to send this to zero as $m \rightarrow \infty$ we use some uniform boundedness of the term $\sum_{i=1}^\infty\norm{\mathcal{G}_i(s,\sy^m_s)}_{H}^2$ which in the case of a Lipschitz operator as in the original paper is immediate from (\ref{galerkinboundstoppingtime}). Where $\mathcal{G}_i$ is a differential operator we must obtain uniform boundedness of the solutions $(\sy^n)$ in a higher norm, hence the need for our space $V$ (which in the context of our SALT Navier-Stokes, would then be $W^{3,2}_{\sigma}(\mathbb{T}^3;\R^3)$). For this reason we must introduce another step to the proof, whereby we show that there exists constants $C,\tilde{C}$ dependent on $M,M',t$ but independent of $n$ such that for the local strong solution $(\sy^n, \tau^{M,t}_n)$ of (\ref{nthorderGalerkin}), \begin{equation} \label{firstresultofuniformbounds}
    \mathbb{E}\sup_{r\in [0,\tau^{M,t}_n]}\norm{\sy^n_{r}}_H^2 + \mathbb{E}\int_{0}^{\tau^{M,t}_n}\norm{\sy^n_s}_V^2ds\leq C\left[\mathbb{E}\left(\norm{\sy^n_{0}}_H^2\right) + 1\right]
\end{equation}
and in particular
\begin{equation} \label{secondresultofuniformbounds}
    \mathbb{E}\sup_{r\in [0,\tau^{M,t}_n]}\norm{\sy^n_{r}}_H^2 + \mathbb{E}\int_{0}^{\tau^{M,t}_n}\norm{\sy^n_s}_V^2ds\leq \tilde{C}.
\end{equation} 
This result is proven by considering $V_n$ as a Hilbert Space with $H$ inner product, applying the It\^{o} Formula in this context and using Assumption \ref{assumptions for uniform bounds2}. (\ref{secondresultofuniformbounds}) then follows from (\ref{firstresultofuniformbounds}) due to (\ref{projectionsboundedonH}) so we see the significance of starting from an initial condition bounded in $H$ and not just $U$ (or at least, square integrable in $H$). From Assumption \ref{therealcauchy assumptions}, along with the requirement that each $\mathcal{P}_n$ is an orthogonal projection in $X$ and $U$ and the conditions (\ref{bilinear form}),(\ref{mu1}),(\ref{mu2}), we deduce that for any $m<n$ and $\lambda_m:= \min\{\mu_m,\mu_m^2\}$,
\begin{align}
\nonumber &2\inner{\mathcal{P}_n\mathcal{A}(t,\boldsymbol{\phi}) - \mathcal{P}_m\mathcal{A}(t,\boldsymbol{\psi})}{\boldsymbol{\phi} - \boldsymbol{\psi}}_U + \sum_{i=1}^\infty\norm{\mathcal{P}_n\mathcal{G}_i(t,\boldsymbol{\phi}) - \mathcal{P}_m\mathcal{G}_i(t,\boldsymbol{\psi})}_U^2\\ &\leq c_{t}\tilde{K}_2(\boldsymbol{\phi},\boldsymbol{\psi}) \norm{\boldsymbol{\phi}-\boldsymbol{\psi}}_U^2 - \frac{\kappa}{2}\norm{\boldsymbol{\phi}-\boldsymbol{\psi}}_H^2 + \frac{c_t}{\lambda_m}K(\boldsymbol{\phi},\boldsymbol{\psi})\left[1 + \norm{\boldsymbol{\phi}}_V^2 + \norm{\boldsymbol{\psi}}_V^2\right] ,\nonumber\\
\nonumber
   &\sum_{i=1}^\infty \inner{\mathcal{P}_n\mathcal{G}_i(t,\boldsymbol{\phi}) - \mathcal{P}_m\mathcal{G}_i(t,\boldsymbol{\psi})}{\boldsymbol{\phi}-\boldsymbol{\psi}}^2_U\\ & \leq c_{t} \tilde{K}_2(\boldsymbol{\phi},\boldsymbol{\psi}) \norm{\boldsymbol{\phi}-\boldsymbol{\psi}}_U^4  + \frac{c_{t}}{\lambda_m} K(\boldsymbol{\phi},\boldsymbol{\psi})\left[1 + \norm{\boldsymbol{\psi}}_V^2\right]. \nonumber
\end{align}
Along with (\ref{secondresultofuniformbounds}) these bounds allow us to conclude that \begin{align}\nonumber &\lim_{m \rightarrow \infty}\sup_{n \geq m}\bigg[\mathbb{E}\sup_{r \in [0,\tau^{M,t}_m \wedge \tau^{M,t}_n]}\norm{\sy^n_{r} - \sy^m_{r}}_U^2 \\& \qquad \qquad \qquad \qquad \qquad  + \mathbb{E}\int_{0}^{\tau^{M,t}_m \wedge \tau^{M,t}_n}\norm{\sy^n_s-\sy^m_s}_H^2ds\bigg] = 0\label{the cauchy prop}\end{align}
again via an application of the It\^{o} Formula for $V_n$ considered as a Hilbert Space with $U$ inner product, on the difference process $\sy^n-\sy^m$. With similar ideas and the Assumption \ref{assumption for prob in V}, we infer that \begin{align} \nonumber &\lim_{S \rightarrow 0}\sup_{n \in \mathbb{N}}\mathbb{P}\bigg(\bigg\{ \sup_{r \in [0,\tau^{M,t}_n \wedge S]}\norm{\sy^{n}_{r}}^2_U \\& \qquad \qquad \qquad \qquad \qquad + \int_0^{\tau^{M,t}_n \wedge S}\norm{\sy^{n}_{r}}^2_Hdr \geq M-1+\norm{\sy^n_0}_{U}^2 \bigg\}\bigg) = 0.\label{prob condition 1} \end{align}
We then apply Theorem \ref{greenlemma} for $\mathcal{H}_1 = H$, $\mathcal{H}_2 = U$ and claim that the resulting pair $(\sy,\tau^{M,t}_{\infty})$ satisfies the additional properties that:
\begin{itemize}
\item $\sy_{\cdot}\mathbf{1}_{\cdot \leq \tau^{M,t}_{\infty}}$ is progressively measurable in $V$;
    \item For $\mathbb{P}-a.e.$ $\omega$, $\sy(\omega) \in C\left([0,T];H\right)$ and $\sy_{\cdot}(\omega)\mathbf{1}_{\cdot \leq \tau^{M,t}_{\infty}(\omega)} \in  L^2\left([0,T];V\right)$ for all $T>0$;
\item $\sy^{n_l}\rightarrow \sy$ holds in the sense that \begin{equation}\label{expectation convergence}\mathbb{E}\left[\sup_{r \in [0,\tau^{M,t}_{\infty}]}\norm{\sy^{n_l}_{r}-\sy_r }^2_U + \int_0^{\tau^{M,t}_{\infty}}\norm{\sy^{n_l}_{r} - \sy_r  }^2_Hdr\right] \longrightarrow 0.\end{equation}
\end{itemize}
Indeed the first two are true from using the uniform boundedness (\ref{secondresultofuniformbounds}) and taking weakly convergent subsequences in the appropriate spaces, then using uniqueness of limits in the weak topology and the embeddings $V \hookrightarrow H \hookrightarrow U$ to identify this limit with $\sy$. The weak convergence preserves the measurability and so the progressive measurability of each $\sy^n$ (from the continuity and adaptedness in $V_n$) is what gives the result here. The final item is then a simple application of the dominated convergence theorem. To conclude that $(\sy,\tau^{M,t}_{\infty})$ is a $V-$valued local strong solution it only remains to show the identity (\ref{identityindefinitionoflocalsolution}), which is done by taking limits of the corresponding terms in (\ref{identityindefinitionoflocalgalerkinsolution}) and applying (\ref{222}), (\ref{333}) alongside the already used assumptions on the $(\mathcal{P}_n)$. We take the limit in $X$ and argue that the identity being satisfied in $X$ is sufficient to conclude the satisfaction of the identity in $U$, given that all integrals can be constructed in $U$ from the regularity of the solution.\\

We have now shown the existence of a $V-$valued local strong solution but for the bounded initial condition (\ref{boundedinitialcondition}). We then show a uniqueness result for such solutions, which is: suppose that $(\sy^1,\tau_1)$ and $(\sy^2,\tau_2)$ are two $V-$valued local strong solutions of the equation (\ref{thespde}) for a given initial condition $\sy_0$. Then for all $s \in [0,\infty)$, \begin{equation} \nonumber\mathbb{P}\left(\left\{\omega \in \Omega: \sy^1_{s \wedge \tau_1(\omega) \wedge \tau_2(\omega) }(\omega) =  \sy^2_{s \wedge \tau_1(\omega) \wedge \tau_2(\omega)}(\omega)  \right\} \right) = 1. \end{equation}
This is proven through applying Assumption \ref{therealcauchy assumptions} in the context of an It\^{o} Formula in $U$ of the difference process of any two solutions. With this uniqueness in place we then conclude the results of Theorem \ref{theorem1} but still for the bounded initial condition, via similar arguments to those used in \cite{holtzziane}. To pass to a general initial condition we consider a sequence of such maximal strong solutions $(\sy^k,\Theta^k)$ corresponding to the bounded initial conditions $(\sy_0\mathbf{1}_{k \leq \norm{\sy_0}_H \leq {k+1}})$ and use the maximality on these pieces to show that the pair $(\sy,\Theta)$ defined at each time $t \in [0,T]$ and $\omega \in \Omega$ by $$ \sy_t(\omega):= \sum_{k=1}^\infty \sy^k_t(\omega)\mathbf{1}_{k \leq \norm{\sy_0(\omega)}_H < k+1}, \quad \Theta(\omega):=\sum_{k=1}^\infty \Theta^k(\omega)\mathbf{1}_{k \leq \norm{\sy_0(\omega)}_H < k+1}$$
is our desired solution for the initial condition $\sy_0$ (where the limit for $\sy$ is in reality just a finite sum). It is clear that for any $\omega$, there exists a $k$ such that $(\sy(\omega),\Theta(\omega))= (\sy^k(\omega),\Theta^k(\omega))$ so the property (\ref{blowupproperty}) follows from the same property in the case of the bounded initial condition. This rounds off our discussion for the proof of Theorem \ref{theorem1}.\end{proof}

In the case where Assumption Set 2 holds, we then look to use the $V-$valued local strong solutions to obtain an $H-$valued local strong solution but now just for a $U-$valued initial condition. At this juncture it is well worth addressing the question of why we consider these distinct types of solution; that is if we wanted an $H-$valued local strong solution then why not restate Assumption Set 1 for the spaces $V$ as $H$, $H$ as $U$ and $U$ as $X$? The reason lies in the application to our stochastic Navier-Stokes Equation, which would then not satisfy the required assumption. This will be discussed more explicitly in Section \ref{section salt navier stokes in abstract framework}.\\

\begin{proof}[Theorem \ref{theorem2}:]
The idea now is to apply this existence result to the sequence of initial conditions $(\mathcal{P}_n\sy_0)$, and apply the same Theorem \ref{greenlemma} argument to the corresponding sequence of solutions. From here we now need to suppose that Assumption Set 2 holds in addition to Asusmption Set 1. In the same manner we start again from a bounded $\sy_0$, this time such that \begin{equation} \label{bounded initial condition in U}
    \norm{\sy_0(\omega)}_U^2 \leq \tilde{M}.
\end{equation}
We could immediately apply Theorem \ref{theorem1} for each initial condition $\mathcal{P}_n\sy_0$, though we want to apply Theorem \ref{greenlemma} for the same spaces $\mathcal{H}_1 = H$ and $\mathcal{H}_2 = U$. Recall that we could not do this immediately for a $U-$valued initial condition and the sequence of Galerkin solutions due to gaining a suitable control on the noise term arising from the difference of the projections. In the present scenario we consider solutions to the unprojected (\ref{thespde}) and so we are not burdened with this difficulty. An application of Theorem \ref{greenlemma} would rely on us being able to conclude that each maximal solution $(\sy^n,\Theta^n)$ corresponding to the initial condition $\mathcal{P}_n\sy_0$ exists up until the stopping time (\ref{tauMtn}) (where the $\sy^n$ notation has now shifted to the above). This is not immediate from Theorem \ref{theorem1}, though we can use similar maximality arguments to extend these solutions to $\tau^{M,t}_n$ at the cost of some regularity. Indeed for these extended solutions we have only the regularity of the $H-$valued solution but with the additional benefit that $\sy_t(\omega)\mathbf{1}_{\cdot \leq \tau^{M,t}_n(\omega)} \in V$ almost everywhere on the product space $\Omega \times [0,\infty)$. This facilitates the use of Assumption \ref{therealcauchy assumptions} in order to show the Cauchy property (\ref{the cauchy prop}), but only via first using an It\^{o} Formula with the bilinear form $\inner{\cdot}{\cdot}_{X \times H}$. We must make this step as the identity for these extended solutions is only satisfied in $X$ hence we cannot use the $U$ inner product. The stochastic integral though can be constructed in $U$ following from Remark \ref{remark analogy}, and the regularity $\sy_t(\omega)\mathbf{1}_{\cdot \leq \tau^{M,t}_n(\omega)} \in V$ allows us to call upon the property (\ref{bilinear form}) so that we can apply Assumption \ref{therealcauchy assumptions}. Without the uniform boundedness (\ref{secondresultofuniformbounds}) for these solutions we need Assumption \ref{assumption for probability in H} instead of just \ref{assumption for prob in V} to deduce (\ref{prob condition 1}). The conclusion of the proof of Theorem \ref{theorem2} then follows identically to that of \ref{theorem1}, now using Assumption \ref{uniqueness for H valued} for the uniqueness part and (\ref{222*}) to show the convergence of the time integral term when justifying that the limiting pair $(\sy,\tau^{M,t}_{\infty})$ obtained from Theorem \ref{greenlemma} is an $H-$valued local strong solution.
\end{proof}

\subsection{SALT Navier-Stokes in the Abstract Framework} \label{section salt navier stokes in abstract framework}

We now briefly comment on the application of this abstract framework to the equation (\ref{converted equation}) in order to conclude the paper. In the previous subsection we have already established the identification of the spaces $$V:= W^{3,2}_{\sigma}(\mathbb{T}^3;\R^3), \quad H:= W^{2,2}_{\sigma}(\mathbb{T}^3;\R^3), \quad U:=W^{1,2}_{\sigma}(\mathbb{T}^3;\R^3), \quad X:= L^2_{\sigma}(\mathbb{T}^3;\R^3)$$
at which point we address the question posed in that subsection as to why we need to make this effort with first the $V-$valued solutions before showing the existence for the $H-$valued ones. That is, why would Assumption Set 1 not hold if we were to shift the spaces from $V$ to $H$, $H$ to $U$ and $U$ to $X$ (with some modifications of the reference to $X$ in Assumption Set 1)? One clear answer is in the treatment of the nonlinear term for (\ref{uniformboundsassumpt1}): for $H=W^{2,2}_{\sigma}(\mathbb{T}^3;\R^3)$ we have the algebra property of the Sobolev Space which affords us a bound $$\norm{\mathcal{L}_{\phi^n}\phi^n}_2 \leq c\norm{\phi^n}_2\norm{\phi^n}_3$$ using the equivalence of the $\norm{\cdot}_2$ and the standard $W^{2,2}$ one. In the $W^{1,2}$ norm we do not have the same luxury and so this nonlinear term cannot be bounded just in terms of the $W^{1,2}$ and $W^{2,2}$ norms as would be required. It is worth noting the significance of using the $\inner{\cdot}{\cdot}_2$ inner product here, as in the same assumption this facilitates the `integration by parts' property for the Stokes Operator in order to gain the additional control we require (i.e. the $-\kappa\norm{\phi^n}_V^2$ term). There is some additional care required then to control the noise terms in these inner products, but this is facilitated by using the same standard cancellation argument that \begin{equation}\label{cancellation}\inner{\mathcal{L}_{\xi_i}\phi}{\phi}_{L^2}=0\end{equation} for $\phi \in W^{1,2}(\mathbb{T}^3;\R^3)$, as well as appreciating that the commutator $[\Delta,B_i]$ is of second order and commuting through the $B_i$ with $\Delta$ until we reduce to a term of the form (\ref{cancellation}). The control (\ref{xi bound}) allows the $\xi_i$ to be effectively ignored in many of these computations, by just pulling them out with the supremum. We refer once more to \cite{Goodair2} for the complete details. Of course it is Theorem \ref{theorem2} which is what translates into our main Theorem of the paper (\ref{NS theorem}), though it is also worth noting that having showed Theorem \ref{theorem1} in this context then we can also say something about the retained regularity of our solutions coming from a more regular initial condition. To really make this point we'd have to say that the maximal times for the different notions of solution were in fact the same, and this is to be addressed in \cite{Goodair2}.

\section{Appendix}

Here we state \cite{holtzziane} Lemma 5.1.

\begin{theorem}\label{greenlemma}
Let $\mathcal{H}_1 \subset \mathcal{H}_2$ be Hilbert Spaces with continuous embedding, and $(\sy^n)$ be a sequence of processes such that for $\mathbb{P}-a.e.$ $\omega$, $\sy^n(\omega) \in C\left([0,T];\mathcal{H}_2\right) \cap L^2\left([0,T];\mathcal{H}_1\right)$ which is a Banach Space with norm $$\norm{\boldsymbol{\psi}}_{X(T)}:= \left(\sup_{r \in [0,T]}\norm{\boldsymbol{\psi}_{r}}^2_{\mathcal{H}_2} + \int_0^T\norm{\boldsymbol{\psi}_{r}}^2_{\mathcal{H}_1}dr\right)^{\frac{1}{2}}.$$ For some fixed $M > 1$ and $t > 0$ define the stopping times $$\tau^{M,t}_n(\omega) := t \wedge \inf\left\{s \geq 0: \norm{\sy^n(\omega)}_{X(s)}^2 \geq M + \norm{\sy^n_0(\omega)}_{\mathcal{H}_2}^2 \right\}$$ and suppose that \begin{equation} \label{required cauchy}\lim_{m \rightarrow \infty}\sup_{n \geq m}\mathbb{E}\norm{\sy^n-\sy^m}_{X(\tau^{M,t}_m \wedge \tau^{M,t}_n)}^2 = 0\end{equation} and $$\lim_{S \rightarrow 0}\sup_{n \in \mathbb{N}}\mathbb{P}\left(\left\{ \norm{\sy^n}_{X(\tau^{M,t}_n \wedge S)}^2 \geq M-1+\norm{\sy^n_0}_{\mathcal{H}_2}^2 \right\}\right) = 0.$$
Then there exists a stopping time $\tau^{M,t}_{\infty}$, a subsequence $(\sy^{n_l})$ and process $\sy = \sy_{
\cdot \wedge \tau^{M,t}_{\infty}}$ such that:

\begin{itemize}
    \item $\mathbb{P}\left(\left\{ 0 < \tau^{M,t}_{\infty} \leq \tau^{M,t}_{n_l}\right)\right\} = 1$;
    \item For $\mathbb{P}-a.e.$ $\omega$, $\sy(\omega) \in C\left([0,\tau^{M,t}_{\infty}(\omega)];\mathcal{H}_2\right) \cap L^2\left([0,\tau^{M,t}_{\infty}(\omega)];\mathcal{H}_1\right)$;
    \item For $\mathbb{P}-a.e.$ $\omega$, $\sy^{n_l}(\omega) \rightarrow \sy(\omega)$ in \\$\left( C\left([0,\tau^{M,t}_{\infty}(\omega)];\mathcal{H}_2\right) \cap L^2\left([0,\tau^{M,t}_{\infty}(\omega)];\mathcal{H}_1\right), \norm{\cdot}_{X(\tau^{M,t}_{\infty}(\omega))} \right).$
\end{itemize}
\end{theorem}

%
%


\begin{thebibliography}{6}
%

\bibitem{general}
Alonso-Orán, D. and Bethencourt de León, A., 2020. On the well-posedness of stochastic Boussinesq equations with transport noise. Journal of Nonlinear Science, 30(1), pp.175-224.

\bibitem{Primitive}
Brzeźniak, Z. and Slavik, J., 2021. Well-posedness of the 3D stochastic primitive equations with multiplicative and transport noise. Journal of Differential Equations, 296, pp.617-676.

\bibitem{cotter1}
Cotter, C., Crisan, D., Holm, D., Pan, W. and Shevchenko, I., 2020. Modelling uncertainty using stochastic transport noise in a 2-layer quasi-geostrophic model. Foundations of Data Science, 2(2), p.173.

\bibitem{cotter2}
Cotter, C., Crisan, D., Holm, D.D., Pan, W. and Shevchenko, I., 2019. Numerically modeling stochastic Lie transport in fluid dynamics. Multiscale Modeling and Simulation, 17(1), pp.192-232.

\bibitem{crisan3deuler}
Crisan, D., Flandoli, F. and Holm, D.D., 2019. Solution properties of a 3D stochastic Euler fluid equation. Journal of Nonlinear Science, 29(3), pp.813-870.

\bibitem{Goodair2}
Crisan, D., Goodair, D. 2022. Analytical Properties of a 3D Stochastic Navier-Stokes Equation. In preparation.

\bibitem{Goodair1}
Crisan, D., Goodair, D., Lang, O., Mensah, P.R., 2022. Existence and Uniqueness of Maximal Strong Solutions to Nonlinear SPDEs with Applications to Viscous Fluid Models. In preparation.

\bibitem{crisanTQG}
Crisan, D., Holm, D.D., Luesink, E., Mensah, P.R. and Pan, W., 2021. Theoretical and computational analysis of the thermal quasi-geostrophic model. arXiv preprint arXiv:2106.14850.

\bibitem{crisanrotating}
Crisan, D. and Lang, O., 2021. Well-posedness Properties for a Stochastic Rotating Shallow Water Model. arXiv preprint arXiv:2107.06601.

\bibitem{crisanlang2}
Crisan, D. and Lang, O., 2021. Local Well-Posedness for the Great Lake Equation with Transport Noise. REV. ROUMAINE MATH. PURES APPL, 66(1), pp.131-155.

\bibitem{crisanlang1}
Crisan, D. and Lang, O., 2022. Well-posedness for a stochastic 2D Euler equation with transport noise. Stochastics and Partial Differential Equations: Analysis and Computations, pp.1-48.

\bibitem {streetvariational}
Crisan, D., and Street, O.D. 2021. Semi-martingale driven variational principles. Proceedings of the Royal Society A, 477(2247), p.20200957.

\bibitem{temamabstract}
Debussche, A., Glatt-Holtz, N. and Temam, R., 2011. Local martingale and pathwise solutions for an abstract fluids model. Physica D: Nonlinear Phenomena, 240(14-15), pp.1123-1144.

\bibitem{holtzziane}
Glatt-Holtz, N. and Ziane, M., 2009. Strong pathwise solutions of the stochastic Navier-Stokes system. Advances in Differential Equations, 14(5/6), pp.567-600.

\bibitem{Goodair stoch}
Goodair, D., 2022. Stochastic Calculus in Infinite Dimensions and SPDEs. arXiv preprint arXiv:2203.17206.

\bibitem {holm2015variational}
Holm, D.D., 2015. Variational principles for stochastic fluid dynamics. Proceedings of the Royal Society A: Mathematical, Physical and Engineering Sciences, 471(2176), p.20140963.

\bibitem {kato1984nonlinear}
Kato, T. and Lai, C.Y., 1984. Nonlinear evolution equations and the Euler flow. Journal of functional analysis, 56(1), pp.15-28.

\bibitem{meminLU}
Mémin, E., 2014. Fluid flow dynamics under location uncertainty. Geophysical and Astrophysical Fluid Dynamics, 108(2), pp.119-146.

\bibitem{robinson}
Robinson, J.C., Rodrigo, J.L. and Sadowski, W., 2016. The three-dimensional Navier–Stokes equations: Classical theory (Vol. 157). Cambridge university press.







\end{thebibliography}
\end{document}